\documentclass[11pt]{article}

\usepackage{amsfonts,amsmath}

\newtheorem {remark} {{\bf Remark}} [section]

\newcommand{\R}{{\mathbb R}}

\def\appendix{\par \setcounter{section}{0} \setcounter{subsection}{0} \def\thesection{\Alph{section}.}}

 \newtheorem{theorem}{Theorem}
 \catcode`@=11
 \@addtoreset{equation}{section}
 
 \catcode`@=12

\begin{document}

\title{A Generalized Occupation Time Formula \\ For Continuous Semimartingales}

\author{{\Large Raouf Ghomrasni} \\School of Computational \& Applied Mathematics,\\University of the Witwatersrand,\\Private Bag 3, Wits, 2050 Johannesburg, South Africa.\\
E-mail: rghomrasni@cam.wits.ac.za}

\date{}

\maketitle

\begin{abstract}
We show that for a wide class of functions $F$ that:
$$
{\lim_{\varepsilon \downarrow 0} \, {\frac{1}{\varepsilon}} \int_0^t
\Big\{ F(s, X_s) \, - \, F(s, X_s - \varepsilon)\Big\} \, d\big<X,X\big>_s} =  -  \int_0^t\!\int_{\R} F(s, x) \, d\, L_s^x \, 
$$
where $X_t$ is a continuous semi-martingale, $(L_t^x, x \in \R, t \geq 0)$ its local time process and $(\big<X,X\big>_t, t \geq 0)$ its quadratic variation process.
\end{abstract}



\section{Introduction}

Recently Feng and Zhao {\cite{FZ}} define  the integral of local time $\int_0^t\!\int_{\R} g(s, x) \, d\, L_s^x$
pathwise and then they derived a generalized It\^o's formula when $\nabla^- F(s,x)$ is only of bounded $p, q$-variation in $(s, x)$. In the case that $g(s, x) = \nabla^- F(s,x)$ is of locally bounded variation in $(s, x)$, the integral $\int_0^t\!\int_{\R} \nabla^- F(s,x) \, d\, L_s^x$ is the Lebesgue-Stieltjes integral. When $g(s, x) = \nabla^- F(s,x)$ is of only locally $p, q$-variation, where $p \ge 1,\, q \ge 1$, and $2q + 1 > 2pq$, the integral is a two-parameter rough path integral rather than a Lebesgue-Stieltjes integral.

In section 2, we first study the time-independent case and establish a formula wich in particular unify the expression of local time of a continuous semimartingales defined as 
\begin{equation} 
L_t^a ={\lim_{\varepsilon \downarrow 0} \, {\frac{1}{\varepsilon}} \int_0^t
1_{[a, a+\varepsilon[}(X_s) \, d\big<X,X\big>_s }  
\end{equation}
and the expression of the quadratic variation process in terms of contributions coming from fluctuations in the process that occur in the vicinity of different spatial points $a \in (-\infty, \infty)$:
\begin{equation} 
 \big<X,X\big>_t =  \int_{\R}  L_t^a \, da
\end{equation}
we then deal with the time-dependent case. A recent survey of semimartingales local time and occupation density concepts is given by I. Serot in {\cite{S}}.

\section{Main Results}

\subsection{Time independent Case}

Using Lyons-Young's integration of one parameter $p$-variation, Feng and Zhao {\cite{FZ}}
defined $\int_{\R} F(x) \, d_x\, L_t^x$ as a rough path integral if $F(x)$ is of
bounded $p$-variation $(1 \leq p < 2)$.
They also proved a dominated convergence theorem
({\cite{FZ}} Theorem 2.1) for the rough path integral and then extended Meyer's formula
to $F(x)$ is of bounded $p$-variation $(1 \leq p < 2)$. We shall use their results in order to establish the following theorem.

\begin{theorem}
Let $F$ be a left continuous function with bounded $p$-variation $(1\leq p < 2)$, we have the following:
\begin{equation} 
{\lim_{\varepsilon \downarrow 0} \, {\frac{1}{\varepsilon}} \int_0^t
\Big\{ F(X_s) \, - \, F(X_s +\varepsilon)\Big\} \, d\big<X,X\big>_s } =  - \int_{\R} F(x) \, d_x\, L_t^x \, 
\end{equation}
\end{theorem}

\begin{remark}
\begin{enumerate}
\item If we take $F(t, x) = 1_{(x \leq a)}$ in (2.1) we have the very definition of $L_t^a$
\item If we take $F(x)=x$ in (2.1) we have $ \big<X,X\big>_t =  \int_{\R}  L_t^x \, dx$.
\end{enumerate}
\end{remark}
{\bf Proof:}
Let us associate to $F$ the following function:
\begin{equation}
H_{\varepsilon}(x):=\frac{1}{\varepsilon} \int_{x}^{x+\varepsilon} F(y)\,dy
\end{equation}
On the one hand we have:
\begin{equation}
H_{\varepsilon}(x):=\frac{1}{\varepsilon} \int_{x}^{x+\varepsilon} F(y)\,dy \rightarrow \quad F(x) \quad \text{for} \quad \varepsilon \to 0
\end{equation}
On the other hand
\begin{equation}
\frac{\partial}{\partial x}H_{\varepsilon}(x):=\frac{1}{\varepsilon} \{F(x+\varepsilon)-F(x)\}
\end{equation}
We note that the function $H_{\varepsilon}(x)$ in (2.2) is of bounded $p$-variation $(1\leq p <2)$ for any fixed $\varepsilon >0$. This is may be easily proved by checking the definition of $p$-variation or as communicated to the author by Prof. Lyons:
\begin{itemize}
\item[(i)] the property of having finite $p$-variation $(p<2)$ can be expressed in
terms of a norm being bounded.
\item[(ii)] the property is preserved under translation.
\item[(iii)] the ball in any norm is convex.
\item[(iv)] the function $H_{\varepsilon}$ is defined as an integral with a convex combination of
translates of the original path.
\end{itemize}
It follows, from Theorem 2.1 in  ({\cite{FZ}})
\begin{equation}
 \int_{\R} H_{\varepsilon}(x) \,d_x\, L_t^x \rightarrow \int_{\R} F(x) \, d_x\, L_t^x \, 
\end{equation}

and 
\begin{equation}
 \int_{\R} H_{\varepsilon}(x) \,d_x\, L_t^x =  \frac{1}{\varepsilon} \int_0^t
\Big\{ F(X_s) \, - \, F(X_s + \varepsilon)\Big\} \, d\big<X,X\big>_s
\end{equation}

We use here Feng-Zhao theorem (2.1) ({\cite{FZ}})

\subsection{Time-dependent Case}

\begin{theorem}
Let $F : [0, t]\times \R \rightarrow  \R$ be a left continuous, locally bounded with bounded $\gamma$-variation in $x$ uniformly in $s$ and of bounded $p, q$-variation in $(s, x)$, where $1 \leq \gamma <2$ and $p, q \geq 1$,
$2q + 1 > 2pq$, Then

\begin{equation} 
{\lim_{\varepsilon \downarrow 0} \, {\frac{1}{\varepsilon}} \int_0^t
\Big\{ F(s, X_s) \, - \, F(s, X_s - \varepsilon)\Big\} \, d\big<X,X\big>_s } =  -  \int_0^t\!\int_{\R} F(s, x) \, d\, L_s^x \, 
\end{equation}

and also,
\begin{equation} 
 {\lim_{\varepsilon \downarrow 0} \, {\frac{1}{2 \, \varepsilon}} \int_0^t
\Big\{ F(s, X_s - \varepsilon) \, - \, F(s, X_s + \varepsilon)\Big\} \, d\big<X,X\big>_s} =  \int_0^t\!\int_{\R} F(s, x) \, d\, L_s^x \,
\end{equation}
\end{theorem}
{\bf Proof:}
By Remark 4.1 and Theorem 4.2 in Feng and Zhao (Two-parameter $p, q$-variation Paths and
Integrations of Local Times):
Let us associate to $F$ the following function:
\begin{equation}
H_{\varepsilon}(t, x):=\frac{1}{\varepsilon} \int_{x}^{x+\varepsilon} F(t,y)\,dy
\end{equation}
On the one hand we have:
\begin{equation}
H_{\varepsilon}(t,x):=\frac{1}{\varepsilon} \int_{x}^{x+\varepsilon} F(t,y)\,dy \rightarrow \quad F(t,x) \quad \text{for} \quad \varepsilon \to 0
\end{equation}
On the other hand
\begin{equation}
\frac{\partial}{\partial x}H_{\varepsilon}(t,x):=\frac{1}{\varepsilon} \{F(t,x+\varepsilon)-F(t,x)\}
\end{equation}
We check easily that the function $H_{\varepsilon}(s,x)$ is of bounded $\gamma$-variation in $x$ uniformly in $s$ and of bounded $p, q$-variation in $(s, x)$, where $1 \leq \gamma <2$ and $p, q \geq 1$,
$2q + 1 > 2pq$ for any fixed $\varepsilon >0$ (see Proof of Theorem 1 above for similar arguments). It follows:

$$
 \int_0^t\!\int_{\R} H_{\varepsilon}(s, x) \,d\, L_s^x \rightarrow  \int_0^t\!\int_{\R} F(s, x) \, d\, L_s^x \, 
$$

and 
$$
 \int_0^t\!\int_{\R} H_{\varepsilon}(s, x) \,d\, L_s^x =  \frac{1}{\varepsilon} \int_0^t
\Big\{ F(s, X_s) \, - \, F(s, X_s + \varepsilon)\Big\} \, d\big<X,X\big>_s
$$

\section{Occupation Time Formula}
 
When $F_x(t,x)=f(t,x )$ exists, (1.1) becomes the  classical occupation time formula for continuous semimartinagles:

\begin{equation} 
 \int_0^t f(s, X_s) \, d\big<X,X\big>_s =  \int_{\R}\!\int_0^t f(s, x) \, d_s\, L_s^x \, dx
\end{equation}

\end{document}